\documentclass{article}

\usepackage{amsmath,amsthm,amssymb}

\newcommand{\mpd}{\;\text{mod}\,p^2}
\newcommand{\mpu}{\;\text{mod}\,p}
\newcommand{\mb}{\mathcal{B}}

\author{Claire Levaillant}
\title{A congruence concerning a convolution involving weighted Bernoulli numbers}
\newcommand{\shnoi}{\sum_{\begin{array}{l}m=1\\\text{m odd}\end{array}}^{p-1}H_m}
\begin{document}

\maketitle

\begin{abstract}
Given a prime $p\geq 5$, we reduce modulo $p$ a convolution of order $p-1$ of powers of two weighted Bernoulli numbers with Bernoulli numbers in terms of harmonic numbers and generalized harmonic numbers. \\Our proof is based on studying the $p$-adic expansion of the number of permutations of $Sym(p-2)$ with an even number of ascents, up to the modulus $p^2$. %We use some congruences due to Emma Lehmer, an identity due to Alzer et al. providing a sum of divided harmonic numbers and several results of \cite{LEV}.
\end{abstract}

\section{Introduction and main result}

% ***Introduction***
Convolutions involving Bernoulli numbers $(B_i)_{i\geq 0}$ or divided Bernoulli numbers $(\mb_i)_{i\geq 1}$ have drawn the attention of mathematicians in the past decades. \\Euler already had come up with an identity regarding a binomial convolution of Bernoulli numbers. This identity reads:
$$\sum_{j=0}^n\binom{n}{j}B_jB_{n-j}=-nB_{n-1}-(n-1)B_n\qquad (n\geq 1)$$
More recently, Miki relates in \cite{MIK} a binomial convolution of divided Bernoulli numbers to a convolution of divided Bernoulli numbers and to harmonic numbers as follows.
$$\sum_{j=2}^{n-2}\binom{n}{j}\mb_j\mb_{n-j}=\sum_{j=2}^n\mb_j\mb_{n-j}-2\mb_n\mathcal{H}_n$$
Since then, many authors have been working on similar congruences, see $\S\,1$ of \cite{DIL}.
Working in the field of $p$-adic numbers, obtaining congruences modulo powers of $p$ for convolutions of divided or ordinary Bernoulli numbers of order $p-k$, with $k$ odd integer, is in general a difficult problem. However, a few special cases corresponding to $k\in\lbrace1,3,5\rbrace$ can be worked out. For instance, it is a straightforward consequence of the Clausen-Von Staudt's theorem and of Euler's identity applied with $n=p-1$ that
$$\sum_{j=2}^{p-3}B_jB_{p-1-j}=1\;\mpu$$
By working with multiple harmonic sums, Zhao has also obtained congruences for convolutions of respective orders $(p-3)$ and $(p-5)$. In \cite{ZHA}, he shows that:
\begin{eqnarray*}
\sum_{j=2}^{p-5}B_jB_{p-3-j}&=&-2B_{p-3}\qquad\qquad\qquad\;\mpu\\
\sum_{j=2}^{p-7}B_jB_{p-5-j}&=&-2B_{p-5}-\frac{2}{3}B_{p-3}^2\qquad\mpu
\end{eqnarray*}
In \cite{LEV3}, we establish similar congruences concerning convolutions of divided Bernoulli numbers of respective orders $p-1$, $p-3$ and $p-5$. Using the notations exposed further below, we show that:
\begin{eqnarray*}
\sum_{j=2}^{p-3}\mb_j\mb_{p-1-j}&=&(2p\mb_{2(p-1)}-p^2\mb_{p-1}^2)_2\qquad\qquad\qquad\qquad\qquad\qquad\;\,\mpu\\
\sum_{j=2}^{p-5}\mb_j\mb_{p-3-j}&=&2((pB_{p-1})_1-1)\mb_{p-3}+2(\mb_{2p-4}-\mb_{p-3})_1\qquad\qquad\,\mpu\\
\sum_{j=2}^{p-7}\mb_j\mb_{p-5-j}&=&-\mb_{p-3}^2+2((pB_{p-1})_1-1)\mb_{p-5}+2(\mb_{2p-6}-\mb_{p-5})_1\mpu
\end{eqnarray*}

\noindent The present paper is concerned with the $p$-residue of a convolution of order $p-1$ of powers of two weighted ordinary Bernoulli numbers with ordinary Bernoulli numbers. Congruences involving harmonic numbers or Bernoulli numbers and powers of two weights have recently appeared in the literature. For instance, congruences for powers of two divided harmonic numbers or generalized harmonic numbers (resp divided Bernoulli numbers) get established in \cite{SUB} (resp in \cite{LEV3}).

\begin{equation*}
\qquad\;\;\sum_{k=1}^{p-1}\frac{H_k}{k2^k}=\frac{7}{24}pB_{p-3}\mpd,\qquad\sum_{k=1}^{p-1}\frac{H_k^{(2)}}{k2^k}=-\frac{3}{8}B_{p-3}\mpu\qquad\qquad\cite{SUB}
\end{equation*}
\begin{equation*}\qquad\qquad\qquad\sum_{k=1}^{p-2}\frac{B_k}{k2^k}=-\frac{1}{2}H_{\frac{p-1}{2}}+\frac{pB_{p-1}+1}{p}-1\mpu \qquad\qquad\qquad\qquad\cite{LEV3}
\end{equation*}

Some key numbers in our proofs are the Eulerian numbers. The Eulerian number $E(n,m)$ is by combinatorial definition the number of permutations on $n$ letters with $m$ ascents. If $(i_1,i_2,\dots,i_n)$ denotes the permutation mapping the integer $j$ onto $i_j$, an ascent is when $i_{k+1}>i_k$.
Eulerian numbers originate in Euler's book dating from $1755$ in which he investigates shifted forms of what are now called Eulerian polynomials:
$$E_n(t)=\sum_{m=0}^n E(n,m)\,t^m$$
Some of the history of the Eulerian numbers appears in \cite{FRO}.
When $n$ is odd, alternating sums of Eulerian numbers are called Euler numbers with odd indices. They relate to the Bernoulli numbers by
$$2^{n+1}(2^{n+1}-1)\frac{B_{n+1}}{n+1}=\sum_{m=0}^{n-1}(-1)^m\,E(n,m)\qquad\qquad (E)_n$$
Taken unsigned (the so-called "Euler tangent numbers"), they count the number of alternating permutations of $Sym(n)$. A permutation is "alternating" if the indices of the descents are all the odd indices. In \cite{LEV}, we show that the Fermat quotient $q_2(p)$ to base $2$ and the Euler number $E_{p-2}$ have the same $p$-residue. Modulo $p$, the Euler number $E_{p-2}$ easily relates to the number of permutations of $Sym(p-2)$ with an even number of ascents, since the non alternating sum of Eulerian numbers equals the order of the symmetric group $Sym(p-2)$ that is $(p-2)!$ and Wilson's theorem asserts that $(p-1)!=-1\mpu$. In \cite{LEV} we denoted the number of permutations of $Sym(k)$ with an even number of ascents by $N_k$. We pushed the expansion of $q_2(p)$ to the modulus $p^2$ in terms of the residue of the Agoh-Giuga quotient and of the first two residues in the $p$-adic expansion of $2N_{p-2}$. We deduced that in order for a prime to be Wilson and super Wieferich ($q_2(p)=0\;\mpd$), it is necessary that the latter two residues are both equal to $1$. Aiming at a search for Wilson and super Wieferich primes, it thus became of interest to investigate the $p$-adic expansion of $N_{p-2}$ to the modulus $p^2$. We then reduced modulo $p^2$ the Eulerian number $E(p-2,2m)$ to a sum of a sum of powers of integers and of a truncated sum of "divided shifted" harmonic numbers. Our novel congruence concerning a convolution of order $p-1$ of powers of two weighted Bernoulli numbers with Bernoulli numbers arises from reducing further the sum over the even ascents of the main two components of $E(p-2,2m)$ just mentioned above. We use several results from \cite{LEV} which we state at the end of this introduction. But first we will introduce a few notations which we shall use throughout the paper. \\

We will start with the harmonic numbers and generalized harmonic numbers. Our notations are the following.
$$H_k:=\sum_{i=1}^k\frac{1}{i}\qquad H_k^{(n)}:=\sum_{i=1}^k\frac{1}{i^n}$$
Second, our notations with respect to sums of powers of integers will be as follows.
$$S_{n,k}:=\sum_{a=1}^n a^k=1^k+2^k+\dots+n^k$$
We now introduce some specific notation with respect to the Bernoulli numbers followed by some notations with respect to $p$-adic numbers. \\
By a theorem due to Clausen \cite{CLA} and independently Von Staudt \cite{VOS}, the denominator of a Bernoulli number $B_k$ consists of products of primes $p$ of multiplicity one, such that $p-1|k$. In particular, $pB_{p-1}=-1\mpu$. In \cite{AGO} (resp \cite{GIU}), Takasaki Agoh (resp Giuseppe Giuga) have conjectured that $nB_{n-1}=-1\;\text{mod}\;n$ if and only if $n$ is a prime. We will denote the Agoh-Giuga quotient by
$$(pB_{p-1})_1:=\frac{1+pB_{p-1}}{p}$$
When $x$ is a $p$-adic integer, we denote by $(x)_i$ its $(i+1)$-th $p$-residue in its Hensel expansion (see e.g. \cite{GOU}), that is:
$$x=\sum_{i=0}^{\infty}(x)_ip^i$$
It is important to underline that this notation does not apply to, nor match with $(pB_{p-1})_1$ which got defined independently by being the Agoh-Giuga quotient. \\
When $a$ is an integer with $1\leq a\leq p-1$, we define a convolution of order $p-1$ of Bernoulli numbers with weighted Bernoulli numbers by:
$$CB_w^{(a)}(p-1):=\sum_{i=2}^{p-3}\frac{B_i}{a^i}B_{p-1-i}$$
The present paper is concerned with $a=2$. \\
Finally, we recall from before that $N_{p-2}$ denotes the number of permutations of $Sym(p-2)$ with an even number of ascents and $q_2(p)$ denotes the Fermat quotient in base $2$. \\

We now state the results of \cite{LEV} in the order in which we will presently use them.
\newtheorem{Result}{Result}

\begin{Result} (Proposition $8$ of \cite{LEV})
$$N_{p-2}=\sum_{m=0}^{\frac{p-3}{2}}\Bigg(S_{2m+1,\,p-2}-p\;\sum_{K=p-(2m+1)}^{p-2}\frac{H_K}{K+(2m+2)}\Bigg)\mpd$$
\end{Result}
\begin{Result} (Theorem $1$ of \cite{LEV})
$$q_2(p)=2N_{p-2}-1\;\mpu$$
\end{Result}
\begin{Result} (Theorem $8(i)$ of \cite{LEV})
$$\sum_{\begin{array}{l}x=1\\\text{x odd}\end{array}}^{p-1} x^{p-2}= (2N_{p-2})_0-1+p\,\Bigg((p\,B_{p-1})_1+(2N_{p-2})_1-\big((2N_{p-2})_0-1\big)^2-2\Bigg)\mpd$$
\end{Result}
\begin{Result} (Proposition $4$ of \cite{LEV})
$$N_{p-2}=\sum_{\begin{array}{l}m=1\\\text{$m$ odd}\end{array}}^{p-2}H_m\;\mpu$$
\end{Result}

\noindent Result $1$ is a refinement of Result $4$ modulo $p^2$. Both results are obtained from reducing a classical formula on Eulerian numbers, see \cite{COM} pp. $243$:
$$E(p-2,2m)=\sum_{k=0}^{2m}(-1)^k\binom{p-1}{k}(2m+1-k)^{p-2}$$
As far as Result $2$, it is obtained from using the definition of $E_{p-2}$ in terms of Bernoulli numbers, some congruences of \cite{LEM} that are due to Emma Lehmer, as well as Eisenstein's formula \cite{EIS} combined with Wolstenholme's theorem \cite{WOL}. \\
Result $3$ gets derived from a congruence originally due to Dmitry Mirimanoff in \cite{MIR} for the sum of odd powers of the first $\frac{p-1}{2}$ integers combined with $(E)_{p-2}$, Result $2$ and a congruence originally due to James Whitbread Lee Glaisher in \cite{GLA} providing $(p-1)!$ modulo $p^2$ as $(p-1)!=pB_{p-1}-p\mpd$.\\

% ***State Main result***
We now state below our main result.

\newtheorem{Theorem}{Theorem}

\begin{Theorem}
Let $H_n$ denote the harmonic number of order $n$ and $H_n^{(k)}$ denote the generalized harmonic number of order $n$. The following congruence holds:
\begin{equation*}
\begin{split}
&\sum_{j=2}^{p-3}\frac{B_j}{2^j}B_{p-1-j}=-1+2\Bigg[\frac{1}{2}\bigg(2\sum_{\begin{array}{l}m=1\\\text{m odd}\end{array}}^{p-1}H_m\bigg)_0\Bigg]_1+\Bigg[2\bigg(\shnoi\bigg)_0\Bigg]_1\\
&+6\shnoi+4\sum_{m=1}^{\frac{p-3}{2}}H_{2m}^{(2)}-4\sum_{m=1}^{\frac{p-3}{2}}H_{2m}H_{2m+1}-4\bigg(\shnoi\bigg)^2\\&
+2\sum_{m=2}^{\frac{p-3}{2}}\Big(\frac{H_1}{2m-1}+\dots+\frac{H_{2m-1}}{1}\Big)\qquad\mpu
\end{split}
\end{equation*}
\end{Theorem}

\newtheorem{Remark}{Remark}
\begin{Remark}
Theorem $1$ is given only in terms of harmonic numbers and generalized harmonic numbers. By Eisenstein's formula for the Fermat quotient
$$q_2(p)=\frac{1}{2}\sum_{k=1}^{p-1}\frac{(-1)^{k-1}}{k}\;\mpu$$
appearing in \cite{EIS} and by Wolstenholme's theorem \cite{WOL} asserting in particular that $H_{p-1}=0\mpd$, we derive that $q_2(p)=H^{'}_{p-1}\mpu$, where $H^{'}_{p-1}$ is the sum of odd reciprocals:
$$H^{'}_{p-1}:=1+\frac{1}{3}+\frac{1}{5}+\dots+\frac{1}{p-2}$$ Then, using Results $2$ and $4$ of the introduction, the sum of harmonic numbers with odd indices may be expressed modulo $p$ in terms of $H^{'}_{p-1}$ as:
$$\sum_{\begin{array}{l}m=1\\\text{$m$ odd}\end{array}}^{p-1}H_m=\frac{H^{'}_{p-1}+1}{2}\mpu$$
\end{Remark}

\section{Proof of the theorem}

We will reduce the main two components of Result $1$. Contrary to the second component of the main sum, for the first component, we will deal directly with the sum over the even ascents. The results are gathered in the following two lemmas. \\Results $2$ and $3$ serve only in the proof for the first lemma. The second lemma has an independent proof which uses in particular a generalization of identities due to Junesang Choi and Hari Mohan Srivastava for sums of "divided shifted" harmonic numbers.\\ Result $4$ will be used only for the final expression like stated in Theorem $1$.
\newtheorem{Lemma}{Lemma}

\begin{Lemma}
\begin{equation*}\begin{split}\sum_{m=0}^{\frac{p-3}{2}}S_{2m+1,p-2}=\frac{1}{2}(2N_{p-2})_0&+p\Big(\frac{1}{2}(2N_{p-2})_0+\frac{1}{2}(2N_{p-2})_1\\&-\frac{1}{2}\big((2N_{p-2})_0-1\big)^2-\frac{1}{2}CB_w^{(2)}(p-1)-1\Big)\mpd
\end{split}\end{equation*}\end{Lemma}
\begin{Lemma}
\begin{equation*}\begin{split}-p\sum_{K=p-(2m+1)}^{p-2}\frac{H_K}{K+(2m+2)}=p&\Big(2H_{2m}^{(2)}-2H_{2m}H_{2m+1}\\&+\frac{H_1}{2m-1}+\frac{H_2}{2m-2}+\dots+\frac{H_{2m-1}}{1}\Big)\mpd
\end{split}\end{equation*}\end{Lemma}

The proof of Lemma $1$ relies on an application of the Bernoulli formula for the sums of powers of integers, see e.g. \cite{KNU}. We have, where we set $B_1=\frac{1}{2}$ instead of $B_1=-\frac{1}{2}$ allowing to avoid the minus signs in the original formula:
$$S_{2m+1,p-2}=\sum_{a=1}^{2m+1}a^{p-2}=\frac{1}{p-1}\sum_{j=0}^{p-2}\binom{p-1}{j}B_j(2m+1)^{p-1-j}\;(B_1=\frac{1}{2})$$
The range of interest for $m$ is such that $1\leq 2m+1\leq p-2$ and so we may write $2m+1=p-2k$ with $k=1,\dots,\frac{p-1}{2}$. Then,
\begin{equation}
\sum_{m=0}^{\frac{p-3}{2}}S_{2m+1,p-2}=-(p+1)\sum_{j=0}^{p-2}\binom{p-1}{j}B_j\sum_{k=1}^{\frac{p-1}{2}}(p-2k)^{p-1-j}\mpd
\end{equation}
We will use several congruences due to Emma Lehmer, each appearing in \cite{LEM}.
\newtheorem*{Theo}{Theorem}
\begin{Theo} (Due to Emma Lehmer \cite{LEM}, $1938$)\\
(i) If $p-1\not|\,2k-2$, then
$$p\,B_{2k}\equiv\frac{1}{2^{2k-1}}\sum_{a=1}^{\frac{p-1}{2}}(p-2a)^{2k}\;\;\text{mod}\,p^3$$
$(ii)$ If  $p-1\not|\,2k-1$, then
$$\sum_{r=1}^{\frac{p-1}{2}}r^{2k}=(2^{-2k+1}-1)B_{2k}\frac{p}{2}\mpd$$
\end{Theo}
\noindent In order to apply point $(i)$ of Lehmer's theorem in Cong. $(1)$ with $2k=p-1-j$, we must have $p-1\not|p-3-j$. This imposes $j\neq p-3$. All the other even values for $j$ satisfy to the condition of application of the theorem. We will write $j=0$, $j=1$ and $j=p-3$ apart. The other values for $j$ provide us with the convolution of interest, up to the upper bound in $j=p-3$, which gets successively added and subtracted.
We obtain:
\begin{equation}\begin{split}
\sum_{m=0}^{\frac{p-3}{2}}S_{2m+1,p-2}=&-(p+1)2^{p-2}pB_{p-1}-\frac{p+1}{2}(p-1)\sum_{k=1}^{\frac{p-1}{2}}(p-2k)^{p-2}\\
&-(p+1)\sum_{j=2}^{p-3}\binom{p-1}{j}B_j2^{p-1-j-1}pB_{p-1-j}\\&+(p+1)\frac{(p-1)(p-2)}{2}B_{p-3}\Big(2pB_2-4\sum_{k=1}^{\frac{p-1}{2}}(k^2-kp)\Big)\;\mpd
\end{split}\end{equation}
After inspection, the corrective term and the term corresponding to $j=p-3$ cancel each other. This is because point $(i)$ of Lehmer's theorem also holds modulo $p^2$ when $2k=2$ but no longer holds modulo $p^3$. Cong. $(2)$ rewrites as:
\begin{equation}\begin{split}
\sum_{m=0}^{\frac{p-3}{2}}S_{2m+1,p-2}=&\frac{1}{2}+\frac{p}{2}-\frac{p}{2}(pB_{p-1})_1+\frac{1}{2}pq_2(p)\\&
+\frac{1}{2}\sum_{k=1}^{\frac{p-1}{2}}(p-2k)^{p-2}-\frac{p}{2}\sum_{j=2}^{p-3}\frac{B_j}{2^j}B_{p-1-j}\;\mpd
\end{split}\end{equation}
Further, we have:
\begin{equation}
\frac{1}{2}\sum_{k=1}^{\frac{p-1}{2}}(p-2k)^{p-2}=-\frac{1}{2}\sum_{\begin{array}{l}k=1\\\text{$k$ even}\end{array}}^{p-1}k^{p-2}-\frac{p}{4}\sum_{k=1}^{\frac{p-1}{2}}k^{p-3}\;\;\;\;\;\mpd
\end{equation}
By the Lehmer congruence $(ii)$ applied with $2k=p-3$, we know that:
$$\sum_{k=1}^{\frac{p-3}{2}}k^{p-3}=0\;\;\mpu$$
Then, the second sum in the right hand side of Cong. $(4)$ vanishes modulo $p^2$. Regarding the first sum, we will use the following result due to Zhi-Hong Sun, combined with our results of $\S\,1$.
\newtheorem*{Lem}{Lemma}
\begin{Lem} (as part of a result by Zhi-Hong Sun in \cite{SU2}, $2000$).\\ Let $k$ be an integer with $k\geq 2$. Then, we have:
$$S_{p-1,k}=p\,B_k+\frac{p^2}{2}\,k\,B_{k-1}\;\mpd$$
\end{Lem}
\noindent Sun's congruence $(5.1)$ of \cite{SU2} for the sums of powers of integers is itself derived from the Bernoulli formula for the sums of powers of integers.\\
The lemma implies that:
$$\sum_{\begin{array}{l}k=1\\\text{$k$ even}\end{array}}^{p-1}k^{p-2}=-\sum_{\begin{array}{l}k=1\\\text{$k$ odd}\end{array}}^{p-1}k^{p-2}\;\;\mpd$$
Then, by applying Results $2$ and $3$ of $\S\,1$, we thus derive Lemma $1$.\\

We now embark onto the proof of Lemma $2$. The proof relies on computing a sum of "shifted divided" harmonic numbers.

First and foremost, an identity for a sum of divided harmonic numbers can be found in \cite{ALZ}. The authors namely prove by induction that:
\begin{equation}
\sum_{j=1}^n\frac{H_j}{j}=\frac{1}{2}\Big[(H_n)^2+H_n^{(2)}\Big]\qquad (n\in\mathbb{N})
\end{equation}

Second, Choi and Srivastava extended this result in \cite{CHO} by showing that:

\begin{eqnarray}
\sum_{j=1}^n\frac{H_j}{j+1}&=&\frac{1}{2}\Big[(H_{n+1})^2-H_{n+1}^{(2)}\Big]\qquad\qquad\qquad\qquad\qquad\qquad\;\; (n\in\mathbb{N})\\
\sum_{j=1}^n\frac{H_j}{j+2}&=&\frac{1}{2}\Big[(H_{n+2})^2-H_{n+2}^{(2)}\Big]+\frac{1}{n+2}-1\qquad\qquad\qquad\;\;\;\; (n\in\mathbb{N})\\
\sum_{j=1}^n\frac{H_j}{j+3}&=&\frac{1}{2}\Big[(H_{n+3})^2-H_{n+3}^{(2)}\Big]+\frac{3}{2(n+3)}+\frac{1}{2(n+2)}-\frac{7}{4}\;\, (n\in\mathbb{N})
\end{eqnarray}

Third, we generalize these identities in the proposition below.

\newtheorem{Proposition}{Proposition}

\begin{Proposition} Let $n$ and $s$ be two integers with $s\geq 3$. Set $H_0:=0$. \\The following identity holds.
\begin{equation*}\begin{split}\sum_{j=1}^n\frac{H_j}{j+s}=&\frac{1}{2}\Big((H_{n+s})^2-H_{n+s}^{(2)}\Big)+\sum_{i=0}^{s-2}\frac{1}{n+s-i}(H_{s-1}-H_i)\\&
-\frac{1}{2}\Big((H_s)^2-H_s^{(2)}\Big)-H_{s-1}H_s+\sum_{k=1}^{s-1}\frac{H_k}{s-k}\end{split}\end{equation*}
\end{Proposition}

\textsc{Proof of Proposition $1$.}
We write:
$$H_{j+s}=H_j+\frac{1}{j+1}+\dots+\frac{1}{j+s}$$
It comes:
\begin{eqnarray}
\sum_{j=1}^{n}\frac{H_j}{j+s}&=&\sum_{j=1}^n\frac{H_{j+s}}{j+s}-\sum_{k=1}^{s-1}\sum_{j=1}^n\frac{1}{(j+k)(j+s)}-\sum_{j=1}^n\frac{1}{(j+s)^2}\\
&=&\sum_{j=s+1}^{n+s}\frac{H_j}{j}-\sum_{k=1}^{s-1}\frac{1}{s-k}\Big(\sum_{j=1}^n\frac{1}{j+k}-\sum_{j=1}^n\frac{1}{j+s}\Big)-\sum_{j=s+1}^{n+s}\frac{1}{j^2}\\
&=&\frac{1}{2}\Big((H_{n+s})^2-H_{n+s}^{(2)}\Big)-\frac{1}{2}\Big((H_s)^2-H_s^{(2)}\Big)-\sum_{k=1}^{s-1}\frac{1}{s-k}\big(H_{n+k}-H_k-H_{n+s}+H_s\big)\notag\\
&&
\end{eqnarray}
where $(11)$ is derived from several applications of $(5)$. \\
Thus,
\begin{equation}\begin{split}\sum_{j=1}^{j+s}\frac{H_j}{j+s}=&\frac{1}{2}\Big((H_{n+s})^2-H_{n+s}^{(2)}\Big)+\sum_{k=1}^{s-1}\frac{1}{s-k}\bigg(\frac{1}{n+k+1}+\dots+\frac{1}{n+s}\bigg)\\
&-\frac{1}{2}\Big((H_s)^2-H_s^{(2)}\Big)-H_sH_{s-1}+\sum_{k=1}^{s-1}\frac{H_k}{s-k}\end{split}
\end{equation}
From there, Proposition $1$ follows. \hfill $\square$\\

We are now in a position to deal with the truncated sum of "shifted divided" harmonic numbers of interest here. First and foremost, we write it as a difference
$$\sum_{j=p-(2m+1)}^{p-2}\frac{H_j}{j+(2m+2)}=\sum_{j=1}^{p-2}\frac{H_j}{j+(2m+2)}-\sum_{j=1}^{p-(2m+2)}\frac{H_j}{j+(2m+2)}$$
We will apply Proposition $1$ with $s=2m+2$ and $2m+2=4,6,\dots,p-1$.
In Proposition $1$, all the terms depending only on $s$ cancel each other within the difference. We thus obtain:
\begin{equation}\begin{split}
\sum_{j=p-(2m+1)}^{p-2}\frac{H_j}{j+(2m+2)}=&\frac{1}{2}\Big(H_{p+2m}^2-H_{p+2m}^{(2)}\Big)-\frac{1}{2}\Big(H_p^2-H_p^{(2)}\Big)\\&+\frac{1}{p+2m}H_{2m+1}+\frac{1}{p+2m-1}(H_{2m+1}-1)\\
&+\frac{1}{p+2m-2}(H_{2m+1}-H_2)+\dots+\frac{1}{p}(H_{2m+1}-H_{2m})\\
&-\frac{1}{p}H_{2m+1}-\frac{1}{p-1}(H_{2m+1}-1)\\
&-\frac{1}{p-2}(H_{2m+1}-H_2)-\dots-\frac{1}{p-2m}(H_{2m+1}-H_{2m})
\end{split}\end{equation}
We get:
\begin{equation}\begin{split}
-p\sum_{j=p-(2m+1)}^{p-2}\frac{H_j}{j+(2m+2)}&=-\frac{p}{2}\Big(H_{p+2m}^2-H_{p+2m}^{(2)}\Big)+\frac{p}{2}\Big(H_p^2-H_p^{(2)}\Big)-2pH_{2m+1}H_{2m}\\&
+p\sum_{j=1}^{2m}\frac{H_j}{j}+p\Big(\frac{H_1}{2m-1}+\frac{H_2}{2m-2}+\dots +\frac{H_{2m-1}}{1}\Big)\\&+H_{2m}\qquad\mpd
\end{split}\end{equation}
By using the Alzer et all identity, Cong. $(14)$ still rewrites as:
\begin{equation}\begin{split}
-p\sum_{j=p-(2m+1)}^{p-2}\frac{H_j}{j+(2m+2)}&=-\frac{p}{2}(H_{p+2m}^2-H_{2m}^2-H_p^2)+\frac{p}{2}(H_{p+2m}^{(2)}+H_{2m}^{(2)}-H_p^{(2)})\\
&-2pH_{2m+1}H_{2m}+p\Big(\frac{H_1}{2m-1}+\frac{H_2}{2m-2}+\dots+\frac{H_{2m-1}}{1}\Big)\\&+H_{2m}\qquad\mpd
\end{split}\end{equation}
%We now compute $H_{p-1-2k}$ modulo $p^2$. We have:
%\begin{eqnarray}
%H_{p-1-2k}&=&1+\frac{1}{2}+\dots+\frac{1}{p-1-2k}\\
%&=&\frac{1}{p-(2k+1)}+\frac{1}{p-(2k+2)}+\dots+\frac{1}{p-(p-1))}\\
%&=&-\frac{1}{(2k+1)^2}(p+(2k+1))-\frac{1}{(2k+2)}(p+(2k+2))-\dots-\frac{1}{(p-1)^2}(p+(p-1))\notag\\
%&&\qquad\qquad\qquad\qquad\qquad\qquad\qquad\qquad\qquad\qquad\qquad\qquad\qquad\qquad\mpd\\
%&=&pH_{2k}^{(2)}+H_{2k}\;\qquad\qquad\qquad\qquad\qquad\qquad\qquad\qquad\qquad\qquad\qquad\mpd
%\end{eqnarray}
%In the congruences above, $(18)$ is obtained from $(17)$ by using Wolstenholme's theorem \cite{WOL} which asserts that:
%$$H_{p-1}=0\mpd\qquad H_{p-1}^{(2)}=0\mpu$$
Moreover, the first factor of the right hand side of $(15)$ may be rewritten as:
$$-\frac{p}{2}\Big\lbrace\bigg(H_{p+2m}+H_{p-1}+\frac{1}{p}\bigg)\bigg(H_{p+2m}-H_p\bigg)-H_{2m}^2\Big\rbrace$$
This in turn reduces modulo $p^2$ as
\begin{equation}-\frac{p}{2}\Bigg\lbrace\Big[2\bigg(H_{p-1}+\frac{1}{p}\bigg)+H_{2m}\Big]\Big[\frac{1}{p+1}+\dots+\frac{1}{p+2m}\Big]-H_{2m}^2\Bigg\rbrace\;\mpd,\end{equation}
By using Wolstenholme's theorem \cite{WOL} which asserts in particular that $$H_{p-1}=0\mpd,$$
$(16)$ reduces further as:
$$-\Big(\frac{1}{p+1}+\dots+\frac{1}{p+2m}\Big)\;\mpd$$
Modulo $p^2$ this sum is: $$pH_{2m}^{(2)}-H_{2m}$$
By gathering all the terms and simplifying, we obtain the result stated in the lemma. \\

By the conjunction of Lemma $1$ and $2$ both plugged into Result $1$ and by using Result $4$, we then derive Theorem $1$.\\\\

\textsc{Email address:} \textit{clairelevaillant@yahoo.fr}

\end{document}